\documentclass[english,12pt]{smfart}
\usepackage{amssymb}
\usepackage{amsfonts}
\usepackage{amscd}
\usepackage[all]{xypic}

\CompileMatrices

\swapnumbers

\def\ord{{\rm ord}}
\def\ac{{\overline{\rm ac}}}

\def\11{{\mathbf 1}}
\def\AA{{\mathbf A}}

\def\CC{{\mathbf C}}

\def\FF{{\mathbf F}}

\def\NN{{\mathbf N}}

\def\PP{{\mathbf P}}
\def\QQ{{\mathbf Q}}
\def\RR{{\mathbf R}}

\def\ZZ{{\mathbf Z}}

\def\cA{{\mathcal A}}
\def\cB{{\mathcal B}}
\def\cC{{\mathcal C}}

\def\cL{{\mathcal L}}
\def\cM{{\mathcal M}}

\def\cO{{\mathcal O}}

\mathchardef\alphag="7C0B \mathchardef\betag="7C0C
\mathchardef\gammag="7C0D \mathchardef\deltag="7C0E
\mathchardef\varepsilong="7C22 \mathchardef\varphig="7C27
\mathchardef\psig="7C20 \mathchardef\zetag="7C10
\mathchardef\epsilong="7C0F \mathchardef\rhog="7C1A
\mathchardef\taug="7C1C \mathchardef\upsilong="7C1D
\mathchardef\iotag="7C13 \mathchardef\thetag="7C12
\mathchardef\pig="7C19 \mathchardef\sigmag="7C1B
\mathchardef\etag="7C11 \mathchardef\omegag="7C21
\mathchardef\kappag="7C14 \mathchardef\lambdag="7C15
\mathchardef\mug="7C16 \mathchardef\xig="7C18
\mathchardef\chig="7C1F \mathchardef\nug="7C17
\mathchardef\varthetag="7C23 \mathchardef\varpig="7C24
\mathchardef\varrhog="7C25 \mathchardef\varsigmag="7C26
\mathchardef\Omegag="7C0A \mathchardef\Thetag="7C02
\mathchardef\Sigmag="7C06 \mathchardef\Deltag="7C01
\mathchardef\Phig="7C08 \mathchardef\Gammag="7C00
\mathchardef\Psig="7C09 \mathchardef\Lambdag="7C03
\mathchardef\Xig="7C04 \mathchardef\Pig="7C05
\mathchardef\Upsilong="7C07

\newtheorem{theorem}[subsection]{Theorem}
\newtheorem{lem}[subsection]{Lemma}
\newtheorem{cor}[subsection]{Corollary}
\newtheorem{prop}[subsection]{Proposition}

\newtheorem{conj}[subsection]{Conjecture}

\theoremstyle{definition}
\newtheorem{definition}[subsection]{Definition}
\newtheorem{example}[subsection]{Example}

\newtheorem{def-prop}[subsection]{Proposition-Definition}
\newtheorem{def-theorem}[subsection]{Theorem-Definition}
\newtheorem{def-lem}[subsection]{Lemma-Definition}

\theoremstyle{remark}
\newtheorem{remark}[subsection]{Remark}

\theoremstyle{plain}

\numberwithin{equation}{subsection}

\def\boxit#1#2{\setbox1=\hbox{\kern#1{#2}\kern#1}%
\dimen1=\ht1 \advance\dimen1 by #1 \dimen2=\dp1 \advance\dimen2 by
#1
\setbox1=\hbox{\vrule height\dimen1 depth\dimen2\box1\vrule}%
\setbox1=\vbox{\hrule\box1\hrule}%
\advance\dimen1 by .4pt \ht1=\dimen1 \advance\dimen2 by .4pt
\dp1=\dimen2 \box1\relax}

\renewcommand{\theequation}{\thesubsection.\arabic{equation}}

\def\AA{{\mathbf A}}

\def\CC{{\mathbf C}}

\def\FF{{\mathbf F}}

\def\NN{{\mathbf N}}

\def\PP{{\mathbf P}}
\def\QQ{{\mathbf Q}}
\def\RR{{\mathbf R}}

\def\ZZ{{\mathbf Z}}

\def\cA{{\mathcal A}}
\def\cB{{\mathcal B}}
\def\cC{{\mathcal C}}

\def\cL{{\mathcal L}}
\def\cM{{\mathcal M}}

\def\cO{{\mathcal O}}

\mathchardef\alphag="7C0B \mathchardef\betag="7C0C
\mathchardef\gammag="7C0D \mathchardef\deltag="7C0E
\mathchardef\varepsilong="7C22 \mathchardef\varphig="7C27
\mathchardef\psig="7C20 \mathchardef\zetag="7C10
\mathchardef\epsilong="7C0F \mathchardef\rhog="7C1A
\mathchardef\taug="7C1C \mathchardef\upsilong="7C1D
\mathchardef\iotag="7C13 \mathchardef\thetag="7C12
\mathchardef\pig="7C19 \mathchardef\sigmag="7C1B
\mathchardef\etag="7C11 \mathchardef\omegag="7C21
\mathchardef\kappag="7C14 \mathchardef\lambdag="7C15
\mathchardef\mug="7C16 \mathchardef\xig="7C18
\mathchardef\chig="7C1F \mathchardef\nug="7C17
\mathchardef\varthetag="7C23 \mathchardef\varpig="7C24
\mathchardef\varrhog="7C25 \mathchardef\varsigmag="7C26
\mathchardef\Omegag="7C0A \mathchardef\Thetag="7C02
\mathchardef\Sigmag="7C06 \mathchardef\Deltag="7C01
\mathchardef\Phig="7C08 \mathchardef\Gammag="7C00
\mathchardef\Psig="7C09 \mathchardef\Lambdag="7C03
\mathchardef\Xig="7C04 \mathchardef\Pig="7C05
\mathchardef\Upsilong="7C07

\DeclareMathOperator*{\grad}{grad}

\def\ord{{\rm ord}}

\newcommand{\abs}[1]{\lvert#1\rvert}

\begin{document}

\title[Igusa's conjecture on exponential
sums]{Igusa's conjecture on exponential sums modulo $p$ and $p^2$
and the motivic oscillation index}

\author[Raf Cluckers]{Raf Cluckers*}
\address{Katholieke Universiteit Leuven,
Departement wiskunde, Celestijnenlaan 200B, B-3001 Leu\-ven,
Bel\-gium. Current address: \'Ecole Normale Sup\'erieure,
D\'epartement de ma\-th\'e\-ma\-ti\-ques et applications, 45 rue
d'Ulm, 75230 Paris Cedex 05, France} \email{cluckers@ens.fr}
\urladdr{www.dma.ens.fr/$\sim$cluckers/}

\subjclass
{Primary 11L07, 
11S40; 
  Secondary 11L05, 
11U09. 
} \keywords{Igusa's conjecture on exponential sums, $p$-adic
exponential sums, finite field exponential sums, Igusa's local zeta
function, $p$-adic cell decomposition, motivic oscillation index,
flaw}

\thanks{*During the realization of this project, the author was a postdoctoral fellow of the Fund for
Scientific Research - Flanders (Belgium) (F.W.O.) and was supported
by The European Commission - Marie Curie European Individual
Fellowship with contract number HPMF CT 2005-007121}

\begin{abstract}
We prove the modulo $p$ and modulo $p^2$ cases of Igusa's conjecture
on exponential sums. This conjecture predicts specific uniform
bounds in the homogeneous polynomial case of exponential sums modulo
$p^m$ when $p$ and $m$ vary. We introduce the motivic oscillation
index of a polynomial $f$ and prove the stronger, analogue bounds
for $m=1,2$ using this index instead of the original bounds.
 The modulo $p^2$ case of our bounds holds for
all polynomials; the modulo $p$ case holds for homogeneous
polynomials and under extra conditions also for nonhomogeneous
polynomials. We obtain natural lower bounds for the motivic
oscillation index by using results of Segers. We also show that, for
$p$ big enough, Igusa's local zeta function has a nontrivial pole
when there are $\FF_p$-rational singular points on $f=0$. We
introduce a new invariant of $f$, the flaw of $f$.
\end{abstract}

\maketitle

\renewcommand{\partname}{}

\section{Introduction}

Let $f$ be a polynomial over $\QQ$ in $n$ variables. Consider the
``global'' exponential sums
$$
E_f(N):= \frac{1}{N^n}\sum_{x\in \{0,\ldots,N-1\}^n}\exp(2\pi
i\frac{ f(x)}{N}),
$$
where $N$ varies over the positive integers. Up to prime
factorization of $N$, to study the dependence of $E_f(N)$ on $N$ it
suffices to study $E_f(p^m)$ in terms of integers $m>0$ and primes
$p$, which requires a mixture of finite field and $p$-adic
techniques.

Igusa's conjecture predicts sharp bounds for $|E_f(p^m)|$ in terms
of integers $m>0$ and primes $p$ when $f$ is homogeneous. Analogue
conjectures for analogue sums over finite field extensions of
$\QQ$ also make sense. Fixing $m=1$ or $m=2$, the conjectured
uniformity in $p$ gives a link between finite field and $p$-adic
exponential sums.

\par
In any case, to find sharp bounds for $|E_f(p^m)|$ for fixed $p$
uniformly in $m$ is often much more easy than bounding uniformly
in $p$ and $m$. Using a resolution of singularities $\pi$ as
below, Igusa defines a rational number $\alpha=\alpha(\pi)\leq 0$
such that
 \begin{equation}\label{a1}
|E_f(p^m)| < c_p m^{n-1} p^{m\alpha}
\end{equation}
for each $m$, each prime $p$, and some constants $c_p$ depending
on $p$. When $f$ is homogeneous, he conjectures in \cite{Igusa3}
that $c_p$ can be taken independently of $p$.

\subsection*{} Let us define the number $\alpha(\pi)$ that is used in
(\ref{a1}) for any polynomial $f$ over $\QQ$, following Igusa
\cite{Igusa3}, or \cite{DenefVeys}, \cite{DenefBour}. If $f$ is a
constant, put $\alpha(\pi):=0$. Otherwise, let
$$\pi_{f}:Y\to \AA^n_\CC$$ be an embedded resolution of
singularities with normal crossings of $f=0$. Let $(N_i,\nu_i)$,
$i\in I$, be the numerical data of $\pi_f$, that is, for each
irreducible component $E_i$ of $f\circ \pi_f=0$, $i\in I$, let
$N_i$ be the multiplicity of $E_i$ in ${\rm div}(f\circ \pi_f)$,
and $\nu_i-1$ the multiplicity of $E_i$ in the divisor associated
to $\pi_f^*(dx)$ with $dx=dx_1\wedge...\wedge dx_n$.
 The \emph{essential numerical data} of $\pi_f$ are
the pairs $(N_i,\nu_i)$ for $i\in J$ with $J=I\setminus I'$ and
where $I'$ is the set of indices $i$ in $I$ such that
$(N_{i},\nu_{i})=(1,1)$ and such that $E_i$ does not intersect
another $E_j$ with $(N_{j},\nu_{j})=(1,1)$.
Define $\alpha(\pi_f)$ as
$$
\alpha(\pi_f)=-\min_{i\in J}\frac{\nu_i}{N_i}
$$
when $J$ is nonempty and define $\alpha(\pi_f)$ as $-2n$ otherwise.
 Take a similar resolution $\pi_{f-c}$ of $f-c=0$ for each critical
value $c\in\CC$ of $f:\CC^n\to \CC$ and put $\pi:=(\pi_{f-c})_c$
where $c$ runs over $\{0\}$ and the critical values of $f$.
 Define
$$
\alpha(\pi):=\max_c (\alpha(\pi_{f-c})),
$$
where $c$ runs over $\{0\}$ and the critical values of $f$.
 Note that the number $\alpha(\pi)$ often depends  on $\pi$, \cite{PrivateVeys}. Namely,
$\alpha(\pi)$ is independent of the choice of $\pi$ when
$\alpha(\pi)\geq - 1$ by properties of the logcanonical threshold of
$f-c$ for all $c$, but on the other hand, when $\alpha(\pi)< -1$,
then $\alpha(\pi)$ always depends on the choice of $\pi$, see Remark
\ref{piVeys}.

\subsection*{}
Now we are ready to state (a slightly generalized form of) what
Igusa conjectured in the introduction of his book of 1978:
\begin{conj}[\cite{Igusa3}]\label{igusacon}
Suppose that $f$  is a homogeneous polynomial over $\QQ$ in $n$
variables. Let $\pi$ and $\alpha(\pi)$ be as above.\footnote{For
homogeneous $f$ one has $\pi =(\pi_f)$ since $0$ is the only
possible critical value of $f$ by Euler's identity.} Then there
exists a constant $c$ such that for each prime $p$ and each positive
integer $m$ one has
$$
|E_{f}(p^m)|  <  c\, m^{n-1} p^{m\alpha(\pi)}.
$$
\end{conj}
This conjecture might hold for a much wider class of polynomials
than homogeneous ones, for example, for quasi-homogeneous
polynomials. It is an open question how to generalize the
conjecture optimally to the case of general $f$, cf.~Example
\ref{exa} and section \ref{sglobal} below.

In certain cases Conjecture \ref{igusacon} is proved. Igusa proved
the case that the projective hypersurface defined by $f$ is smooth
\cite{Igusa3}, by using Deligne's bounds in \cite{DeligneWI} for
finite field exponential sums.
 The case that $\pi$ is a toric resolution and $f$ is
nondegenerate with respect to its Newton polyhedron at the origin
is proven under special conditions by Denef and Sperber in
\cite{DenSper}, and in full generality by the author in
\cite{CDenSper}.

\subsection*{}
In this paper we prove Igusa's Conjecture \ref{igusacon} when we fix
$m=1$, resp.~$m=2$, see Corollary \ref{ci}. Actually, we
prove a bit more:\\

\begin{quote}

\noindent 1) instead of using $\alpha(\pi)$, we define and use an
exponent $\alpha_f$, the \emph{motivic oscillation index of $f$},
which is intrinsically associated to $f$ and at least as sharp and
sometimes sharper
than $\alpha(\pi)$ since $\alpha_f\leq \alpha(\pi)$;\\

\noindent2) for $m=2$ the conjecture holds for arbitrary $f$;
\\

\noindent 3) when we restrict to big enough primes $p$ we prove the
analogue bounds for all finite field extensions of the fields
$\QQ_p$ and all fields $\FF_q((t))$ of characteristic $p$.\\
\end{quote}

The motivic oscillation index $\alpha_f$ is defined and studied in
section \ref{sm}; this notion is based on a suggestion made to us
by Jan Denef. In particular, we show that (\ref{a1}) holds with
$\alpha=\alpha_f$ instead of $\alpha=\alpha(\pi)$.

\subsection*{}
 For arbitrary polynomials $f$ in $n$ variables over some number field,
 we establish in Theorem
\ref{est} the lower bound
 $$
\frac{-n+\delta_f}{2}\leq \alpha_f,
 $$
where $\delta_f$ is the dimension of the locus of $\grad f$, where
we say that the empty scheme has dimension $-\infty$. It is proven
by showing in Theorem \ref{tpol} that for big enough $p$ Igusa's
local zeta function has a nontrivial pole if and only if there is a
$\FF_p$-rational point on $\grad f=f=0$ (Theorem \ref{tpol}) and by
using the lower bounds of Segers \cite{Segers} for the isolated
singularity case in combination with a Cartesian product argument on
integrals.

Having this lower bound for $\alpha_f$, Conjecture \ref{igusacon}
with $m=1$ fixed follows from upper bounds given by Katz
\cite{Katz} for exponential sums over finite fields of large
characteristic. The case $m=2$ follows elementarily from Hensel's
Lemma and bounds on the number
of $\FF_q$-rational points on $\grad f=0$.\\
The analogue results over the fields $\FF_q((t))$ of big enough
characteristic just follow from a comparison principle coming from
the theory of motivic integration, named transfer principle,
cf.~\cite{CLexp}, \cite{CLoesIII}, or with some work from
\cite{DenefBour}.
\\

In Theorem \ref{gig} we generalize our results for homogeneous
polynomials to arbitrary polynomials satisfying Condition
(\ref{cond}) which relates the singluarities at infinity in $\PP^n$
with the singularities in affine space.\\

In section \ref{sglobal} we slightly sharpen Igusa's conjecture by
using $\alpha_f$ instead of $\alpha(\pi)$. In general, $\alpha_f$
might not be the only ingredient for sharp bounds, so we introduce a
new invariant of a polynomial $f$, the flaw of $f$, for the study of
uniform bounds for $|E_f(p^m)|$.

\begin{remark}\label{piVeys}
When $\alpha(\pi)< -1$, then $\alpha(\pi)$ depends on the choice of
$\pi$, as is communicated to the author by Segers and Veys
\cite{PrivateVeys}. Namely, let $E$ be an exceptional irreducible
component (for $\pi$) with numerical data $(N,\nu)$ which intersects
an irreducible component $S$ of the strict transform. The numerical
data of $S$ is automatically of the form $(N_S,1)$, since it is a
component of the strict transform. Since $\alpha(\pi)< -1$, we find
that $N_S=1$ and that $-\nu/N<-1$. Performing one more blowing up
with center the intersection of $E$ and $S$ yields a new exceptional
component with numerical data $(N+1,\nu+1)$, which still intersects
an irreducible component of the strict transform. Repeating this
process with finitely many steps yields an alternative tuple of
resolutions $(\pi')$ with $\alpha(\pi)<\alpha(\pi')<-1$, and with
moreover $\alpha(\pi')$ as close to $-1$ as one wants.

\end{remark}

\subsection*{Some notation}

Let $k$ be a number field and let $f$ be a polynomial over $k$ in
$n$ variables. Define $\alpha(\pi)$ for a tuple $\pi$ of embedded
resolutions of polynomials $f-c$ as in the introduction. Let
$\cA_k$ be the collection of all finite field extensions of
$p$-adic completions of $k$. For $K\in\cA_k$, let $\cO_K$ be its
valuation ring with maximal ideal $\cM_K$ and residue field $k_K$
with $q_K$ elements, and let $\psi_K$ be any additive
character\footnote{For example, one can take the standard
nontrivial additive character
$$\psi_K:K\to\CC:x\mapsto \exp(2\pi i(\mathrm{Trace}_{K,\QQ_p}(x)\bmod
\ZZ_p))
$$
for $K$ a finite field extension of $\QQ_p$, and where
$\mathrm{Trace}_{K,\QQ_p}$ is the trace of $K$ over $\QQ_p$.} from
$K$ to $\CC^\times$ which is trivial on $\cO_K$ and nontrivial on
some element of order $-1$.
 Write $|\cdot|_K$ for the norm on $K$ and
$v_K:K^\times\to\ZZ$ for the valuation.

A Schwartz-Bruhat function $K^n\to \CC$ is a locally constant
function with compact support. For $K$ in $\cA_k$, $\varphi$ a
Schwartz-Bruhat function on $K^n$, and $y\in K$, consider the
exponential integral
$$
E_{f,K,\varphi}(y):= \int_{K^n}\varphi(x)\psi_K(yf(x))|dx|_K,
$$
with $|dx|_K$ the normalized Haar measure on $K^n$, and write
$E_{f,K}$ for $E_{f,K,\varphi_{{\rm triv}}}$ with $\varphi_{{\rm
triv}}$ the characteristic function of $\cO_K^n$.\footnote{For
$k=\QQ$, $K=\QQ_p$, and $\psi_K$ the standard nontrivial additive
character, one thus has $E_{f,K}(p^{-m})=E_f(p^m)$.}

Let $\cB_k$ be the collection of all fields $K$ of the form
$\FF_q((t))$ such that $K$ is an algebra over $\cO_k$.
 Let $\cC_k$ be the union of $\cA_k$ and $\cB_k$.
 Write $\cC_{k,M}$ for
the set of fields in $\cC_k$ with residue field of characteristic
$>M$ and write similarly for $\cA_{k,M}$ and $\cB_{k,M}$.
  For $K$ in $\cB_k$, use the notations $\cO_K$,
$\cM_K$, $\psi_K$, etc., similarly as for $K$ in $\cA_k$. We can
consider $f$ as a polynomial over $\cO_k[1/N]$ for some large enough
integer $N$, and hence, for $K\in\cB_{k,N}$, the polynomial $f$
makes sense as a polynomial over $K$ and $ E_{f,K,\varphi}$ is
defined as the corresponding integral for any Schwartz-Bruhat
function $\varphi$ on $K^n$; for $K\in\cB_k$ of characteristic $\leq
N$, $E_{f,K,\varphi}$ is defined to be zero.

The Vinogradov symbol $\ll$ has its usual meaning, namely that for
complex valued functions $f$ and $g$ with $g$ taking non-negative
real values $f\ll g$ means $|f|\leq c g $ for some constant $c$.

By \cite{DenefVeys}, Proposition 2.7, one has for any $K$ in
$\cA_k\cup\cB_{k,M}$ for big enough $M$ that
 \begin{equation}\label{a2}
|E_{f,K,\varphi}(y)| \ll |v_K(y)|^{n-1} |y|^{\alpha(\pi)}
\end{equation}
for any $\pi$ as above, where both sides are considered as functions
in $y\in K^\times$ and $|v_K(y)|$ denotes the real absolute value of
$v_K(y)$.

Let $C_f$ be the closed subscheme of $\AA^n_k$ given by $0=\grad f$
and let $\delta_f$ be the dimension of $C_f$, where we say that the
empty scheme has dimension $-\infty$. Let $S_f$ be the closed
subscheme of $\AA^n_k$ given by $f=0=\grad f$ and let
$\varepsilon_f$ be the dimension of $S_f$.
 If $f$ is a polynomial over $\cO_k[1/N]$ for some $N$, then for any
ring $R$ which is an algebra over $\cO_k[1/N]$, denote the set of
$R$-rational points on $C_f$ by $C_f(R)$, and likewise on $S_f$. Use
similar notation for polynomials over other fields than $k$.

\section{Igusa's local zeta function has nontrivial
poles}\label{secnont}
Following Shavarevitch, Borevitch, Igusa, and
others, for $f$ a nonconstant polynomial in $n$ variables over
$\cO_k$ and for $K$ in $\cC_k$, consider the Poincar\'e series
$$
P_{f,K}(T):=\sum_{m\in \NN} N_{f,K,m} T^m
$$
with
$$
N_{f,K,m}:=\frac{1}{q_K^{mn}}\sharp \{x\in (\cO_K\bmod
\cM_K^m)^n\mid f(x)\equiv 0 \bmod \cM_K^m\}.
$$

For $K\in\cA_k$ it is known that $P_{f,K}(T)$ is rational in $T$,
see \cite{Igusa3} and \cite{Denef}. Also for $K\in \cB_{k,M}$ for
$M$ big enough, $P_{f,K}(T)$ is rational in $T$, shown using
resolutions with good reduction \cite{Denefdegree},
\cite{DenefBour}, or by insights of motivic integration,
cf.~\cite{DL} or \cite{CLexp}.\footnote{It is conjectured that
$P_{f,K}(T)$ is rational for all $K\in\cB_k$.}

\begin{definition}\label{nont}
For $K$ in $\cC_\QQ$ and $Q_{K}(T)$ a rational function over $\QQ$
in $T$ define the \emph{nontrivial poles} of $Q_{K}(T)$ as the poles
of $Q_{K}(T)$ without the pole $T=q_K$ if this is a pole of
multiplicity exactly $1$.
\end{definition}

\begin{theorem}\label{tpol}
Let $f$ be a nonconstant polynomial in $n$ variables over $k$. There
exists a number $M$ such that for all $K\in\cC_{k,M}$, the rational
function $P_{f,K}(T)$ has a nontrivial pole if and only if
$S_{f}(k_K)$ is nonempty, with $S_f(k_K)$ as in the section on
notation.
\end{theorem}
\begin{proof}

For $M$ big enough and for $K$ in $\cA_{k,M}$, if $S_f(k_K)$ is
nonempty, then $S_f(\cO_K)$ is nonempty and thus $N_{f,K,m}\not =0 $
for all $m$. Hence, $P_{f,K}(T)$ has at least one pole.

In \cite{Denefdegree} it is shown that for $M$ big enough and
$K\in\cA_{k,M}$, the degree of $P_{f,K}(T)$ in $T$ is $\leq 0$.
Hence, the same holds for $K\in \cB_{k,M}$ with $M$ big enough. So,
in the case that $P_{f,K}(T)$ has only $T=q_K$ as pole, if this pole
has multiplicity $1$, and if $K$ lies in $\cC_{k,M}$ with $M$ big
enough, we can write
$$
P_{f,K}(T)= (A_K+B_KT)/(1-q_K^{-1}T),
$$
for some rational numbers $A_K$ and $B_K$ depending on $K$. By
developing the right hand side into series we get
$$
N_{f,K,1}=q_K^{-1}A_K+B_K
$$
and
$$
N_{f,K,2}=q_K^{-2}A_K + q_K^{-1}B_K=q_K^{-1} N_{f,K,1}.
$$
We will only use that $q_K^{-1} N_{f,K,1}=N_{f,K,2}$. Let $\bar f$
be the reduction of $f$ modulo $\cM_K$. Put
$$Y_{K,f}:=\{x\in
k_K^n\mid \bar f(x)=0\not=(\grad \bar f )(x)\}.$$ Suppose now that
$S_f(k_K)$ is nonempty.
 Count points to get
$$
N_{f,K,1}= \frac{1}{q_K^{n}}(\sharp Y_{K,f} + \sharp S_f(k_K))
$$
and
$$
N_{f,K,2}= \frac{1}{q_K^{n+1}}\sharp Y_{K,f}+ \frac{1}{q_K^n}\sharp
S_f(k_K),
$$
by Taylor expansion of $f$ around points with residue in $S_f(k_K)$
resp.~in $Y_{K,f}$, and since the projection
$$
\{x\in (\cO_k/\cM_K^{2}\cO_k)^n\mid f(x)\equiv 0 \bmod \cM_K^{2}\}
 \to
\{x\in (\cO_k/\cM_K\cO_k)^n\mid f(x)\equiv 0 \bmod \cM_K\}
$$
is surjective for $K$ in $\cA_{k,M}$ and $M$ big
enough\footnote{This is so because the projection $\{x\in
\big(k[t]/(t^{2})\big)^n\mid f(x)\equiv 0 \bmod t^{2}\}
 \to
\{x\in k^n\mid f(x)= 0\}$ is surjective for $k$ any field of
characteristic zero. Alternatively, this follows from Greenleaf's
theorem, cf.~\cite{BIRCHMcCANN}.}. This gives a contradiction with
$q_K^{-1} N_{f,K,1}=N_{f,K,2}$.

\end{proof}


\begin{remark}\label{tpz}
By an observation of Igusa, cf.~\cite{DenefBour}, section 1.2, one
has
\begin{equation}\label{PZ}
P_{f,K}(T)=\frac{1-TZ_{f,K}(s)}{1-T},
\end{equation}
 with $T=q_K^{-s}$ and
$Z_{f,K}(s)$ Igusa's local zeta function defined for real $s> 0$ by
$$
Z_{f,K}(s):=\int_{\cO_K^n}|f(x)|^{s}|dx|,
$$
with $|dx|$ the normalized Haar measure on $\cO_K^n$. Since
$Z_{f,K}(0)=1$, $Z_{f,K}(s)$ and $ P_{f,K}(T)$ have exactly the same
poles in $T=q_K^{-s}$ and thus Theorem \ref{tpol} also applies to
$Z_{f,K}(s)$.
\end{remark}

\section{Definition and basic properties of the oscillation
indices}\label{sm}
\subsection{}\label{letfbe} Let $f$ be a polynomial in $n$ variables over a number
field $k$, let $K$ be in $\cC_k$, and let $\varphi:K^n\to \CC$ be a
Schwartz-Bruhat function. Define $\alpha_{f,K,\varphi}$ as the
infimum in $\RR\cup\{-\infty\}$ of all real numbers $\alpha$
satisfying
$$
E_{f,K,\varphi}(y)\ll |y|_K^\alpha,
$$
where both sides are considered as functions in $y\in K^\times$.
Note that $\alpha_{f,K,\varphi}\leq 0$ because $\alpha=0$ always
satisfies the previous Vinogradov inequality. Define
$\alpha_{f,K}$ as the supremum of the $\alpha_{f,K,\varphi}$ when
$\varphi$ varies over all Schwartz-Bruhat functions on $K^n$. Call
$\alpha_{f,K}$ the \emph{$K$-oscillation index of $f$}.

\begin{definition} Define
the \emph{motivic oscillation index of $f$} (over $k$) as
$$
\alpha_f:=\lim\limits_{M\to+\infty}\sup\limits_{K\in\cA_{k,M}}\alpha_{f,K}.
$$
\end{definition}

\begin{lem}\label{lem1}
For $K$ in $\cA_k\cup\cB_{k,M}$ with $M$ big enough and $\varphi$ a
Schwartz-Bruhat function on $K^n$, the number $\alpha_{f,K,\varphi}$
equals $-\infty$ or one of the finitely many quotients $-\nu_i/N_i$
with $(N_i,\nu_i)$ the essential numerical data of the tuple of
resolutions $\pi=(\pi_{f-c})_c$ chosen as in the introduction.
\end{lem}
\begin{proof}
This follows from Proposition 2.7 of \cite{DenefVeys},
cf.~(\ref{Dveys}) below.
\end{proof}
\begin{cor}\label{cor1}
For all $K$ in $\cA_k\cup\cB_{k,M}$ with $M$ big enough there exists
$\varphi$ such that
$$\alpha_{f,K,\varphi}=\alpha_{f,K},$$ and for all $N$ there exists
$K$ in $\cA_{k,N}$ and $L$ in $\cB_{k,N}$ such that
$$\alpha_f=\alpha_{f,K}=\alpha_{f,L}.$$
Moreover, the following list of a priori inequalities holds
\begin{equation}\label{lin}
-\infty\leq \alpha_{f,K,\varphi}\leq \alpha_{f,K}\leq \alpha_{f}\leq
\alpha(\pi)\leq 0
\end{equation}
for big $M$, any $K$ in $\cC_{k,M}$ and any Schwartz-Bruhat function
$\varphi$ on $K^n$.
 \end{cor}
\begin{proof}
The existence of $L$ follows from the motivic understanding of
exponential integrals over $K^n$ for $K$ varying in $\cC_{k,M}$ for
$M$ big enough, as established in \cite{CLexp} (equivalently, one
can use \cite{DenefBour}). The inequality $\alpha_{f}\leq
\alpha(\pi)$ follows from (\ref{a2}).
\end{proof}

Although the $\alpha_{f,K,\varphi}$ are defined as infima, they
yield actual bounds, analogous to (\ref{a1}) and (\ref{a2}):

\begin{lem}\label{lb}
For $K$ in $\cA_k\cup\cB_{k,M}$ for sufficiently big $M$ and
$\varphi:K^n\to \CC$ a Schwartz-Bruhat function, either
$$\alpha_{f,K,\varphi}=-\infty\ \mbox{ and }\
E_{f,K,\varphi}(y)=0\ \mbox{ for  all $y$ with $|y|_K$ sufficiently
big,}
$$
 or,
$$\alpha_{f,K,\varphi}>-\infty\ \mbox{ and }\
E_{f,K,\varphi}(y)\ll |v_K(y)|^{n-1} |y|_K^{\alpha_{f,K,\varphi}}\
\mbox{ as functions in $y\in K^\times$,}
$$
where $|v_K(y)|$ denotes the real absolute value of $v_K(y)$.
\end{lem}
\begin{proof}
The  statements for $K$ in $\cA_k$ are immediate corollaries of the
asymptotic expansions of $E_{f,K,\varphi}(y)$ given by Igusa
\cite{Igusa3} and generalized by Denef and Veys in \cite{DenefVeys},
for $|y|$ going to $\infty$. Namely, by \cite{DenefVeys},
Proposition 2.7, for $K\in\cA_k$ and for $y$ with $|y|$ big enough,
one can write $ E_{f,K,\varphi}(y)$ as a finite $\CC$-linear
combination of terms of the form
\begin{equation}\label{Dveys}
\chi(y)|y|^\lambda v_K(y)^{\beta}\psi(cy)
\end{equation}
with $c\in K$ a critical value of $f$, $\chi$ a character with
finite image of $K^\times$ into the complex unit circle, $\lambda$ a
complex number with negative real part, and $\beta$ an integer
between $0$ and $n-1$, where real parts of the occurring $\lambda$
are among the quotients $-\nu_i/N_i$ with $(N_i,\nu_i)$ the
essential numerical data of some $\pi$ as above. From this the
statements follow for $K\in \cA_k$.

For $K$ in $\cB_{k,M}$ one then uses the transfer principle of
\cite{CLexp}. Namely, by this transfer principle it follows that for
$M$ big enough and for $K$ in $\cB_{k,M}$, $ E_{f,K,\varphi}(y)$ can
be written as a similar linear combination as given higher up in the
proof for $K$ in $\cA_k$. (The statement about $\cB_{k,M}$ also
follows from \cite{DenefBour} and the existence of embedded
resolutions of $f=0$ with good reduction modulo $p$ for $p$ any
prime bigger than some $M$, see also \cite{Denefdegree}.)
\end{proof}

The following Proposition essentially follows from Theorem
\ref{tpol} and \cite{DenefVeys}, with $C_f$ as in the section on
notation.
\begin{prop}[Nontriviality of $\alpha_f$]\label{nontrivialosc}

\item[(i)] For $K$ in $\cA_k$ or in $\cB_{k,M}$ for $M$ big enough,
$\alpha_{f,K}=0$ if and only if $f$ is a constant.

\item[(ii)] If $f$ is not a constant, then there exists $M>0$ such
that, for $K$ in $\cC_{k,M}$,
 $\alpha_{f,K,\varphi_{\rm triv}}=-\infty$ if
and only if $C_{f}(k_K)$ is empty.

By consequence, $-\infty<\alpha_f<0$ if and only if $f$ is
nonconstant and $\grad f=0$ has a solution in $\CC^n$.
\end{prop}
\begin{proof}
The first statement for $K$ in $\cA_k$ follows from (\ref{a2}),
(\ref{lin}) and the definition of $\alpha(\pi)$ in the introduction.
The second statement follows from Theorem \ref{tpol} applied to some
$K$ and one of the polynomials $f-c$ with $c\in \CC$ a critical
value of $f$ with $c\in K$ and from the fact that any nontrivial
pole of $P_{f,K}(T)$ effectively appears as $\lambda$ in one of the
terms (\ref{Dveys}) when writing $ E_{f,K}(y)$ as a linear
combination of terms as (\ref{Dveys}), by Proposition 2.7 of
\cite{DenefVeys}.
\end{proof}

\section{A lemma on parameterized integrals}\label{sc}

Let $\cL_{DP}$ be the language of Denef-Pas, consisting of three
sorts: valued field, residue field and value group, and the language
of fields for the residue field, the language of fields for the
valued field, the language of ordered additive groups for the value
group, and the extra symbols $\ac$ and $\ord$.

Let $k$ be a number field. Consider the language $\cL_{DP}(k)$,
which is the language $\cL_{DP}$ together with coefficients for $k$
in the valued and residue field sort. On a field $K\in\cA_{k}$, the
symbols of $\cL_{DP}(k)$ have their natural meaning, where $\ac(x)$
for nonzero $x\in K$ is interpreted as $x\pi_K^{-\ord x}\bmod \cM_K$
in the residue field $k_K$, for some (consequent) choice of
uniformizer $\pi_K$ of $\cO_K$ and where $\ac(0)=0$. Note that the
language $\cL_{DP}$ can be naturally interpreted in $L((t))$ for any
field $L$ in a similar way.

Then, by quantifier elimination results by Pas \cite{Pas}, for any
given $\cL_{DP}(k)$-formula $\varphi$ there exists a
$\cL_{DP}(k)$-formula $\varphi'$ without valued field quantifiers
such that for all $K\in\cA_{k,M}$ with $M$ big enough, $\varphi$ and
$\varphi'$ determine the same set over $K$.

By a \emph{basic step function} on $K^n$ with $K$ in $\cA_k$ we mean
a characteristic function of a Cartesian product of factors of the
form $\cO_K$ and $a+\cM_K$ for $a\in\cO_K$. For example, the
characteristic function of $\cO_K\times \cM_K\times (1+\cM_K)^{n-2}$
is a basic step function on $K^n$.

The goal of this section is to prove the following lemma.
\begin{lem}\label{padicf}
Let $k$ be a number field and let $f$ be a polynomial in the
variables $x_1,\ldots,x_n$ over $k$. Write $\hat x$ for
$x_2,\ldots,x_n$. Then there are $M>0$ and a polynomial $g$ in one
variable over $k$ such that for all $a\in k$ with $g(a)\not=0$, all
$K\in\cA_{k,M}$, and all basic step functions $\varphi$ on
$K^{n-1}$,  one has that
$$
I_K(b,s,\varphi):=\int_{\cO_K^{n-1}}\varphi(\hat x)|f(b,\hat
x)|^s|d\hat x|,
$$
with $s>0$ a real variable and $|d\hat x|$ the normalized Haar
measure, is independent of $b\in a+\cM_K$.
\end{lem}
\begin{proof}

We give a proof based on motivic integration and constructible
functions as developed in \cite{CLoes}. By \cite{CLoes} (or already
with some work by Pas \cite{Pas}), there are $\cL_{DP}(k)$-formulas
such that their interpretation in $K$ with $K$ a field in
$\cA_{k,M}$ with $M$ big enough yields functions
$$
\alpha_i(K),\beta_i(K):K\to \ZZ
$$
and definable subsets
$$
A_i(K)\subset K\times k_K^{(n-1)+\ell}
$$
of $K$ times a power of the residue field of $K$ such that the
following holds for $K$ in $\cA_{k,M}$ with $M$ big enough: with
$\varphi_z$ the characteristic function of the product of the
residue classes mod $\cM_K$ of $z=(z_2,\ldots,z_n)\in k_K^{n-1}$,
one has that $I_K(b,s,\varphi_z)$ for $b\in K$ is a $\QQ$-linear
combination of products of factors of the form
$$
a_{iz}(K):=\sharp \big(A_i(K)\cap \{b\}\times\{z\}\times
k_K^\ell\big),
$$
$$
q_K^{\alpha_i(K)s+\beta_i(K)},
$$
$$
\frac{1}{1-q_K^{as+b}},
$$
and
$$
\beta_i(K),
$$
with integers $(a,b)\not=(0,0)$. We show that the $a_{iz}(K)$,
$\alpha_i(K)$, and $\beta_i(K)$, are constant on $a+\cM_K$ for all
$a$ with $g(a)\not = 0$ for some polynomial $g$ over $k$, uniformly
in $K$ in $\cA_{k,M}$. But this follows at once from Lemma
\ref{lcd}. For basic step functions which also involve the
characteristic function of Cartesian products with $\cO_K$ one works
similarly, say, with some more (witnessing) $z$-variables. The Lemma
then follows.
 \end{proof}

\begin{definition}[One dimensional Cells]\label{cell1}Let $L$ be any
field of characteristic zero. Let $X\subset L((t))$ be a
$\cL_{DP}(k)$-definable set and
$$
f:X\to S
$$
$$
c:S\to L((t))
$$
$\cL_{DP}(k)$-definable functions, with $S$ a Cartesian product of
some power of $L$ and some power of $\ZZ$.  The tuple $(X,f,c)$ is
called a $(1)$-cell with presentation $f$ and center $c$ when for
each $s\in f(X)$, there exists $\xi\in L^\times$ and $a\in\ZZ$ such
that the fiber $f^{-1}(s)$ is an open ball of the form
$$
\{x\in L((t))\mid \ac(x-c(s))=\xi \wedge \ord (x-c(s))= a\}.
$$
The tuple $(X,f,c)$ is called a $(0)$-cell with presentation $f$ and
center $c$ when $c\circ f$ is the identity on $X$.
\end{definition}

\begin{lem}[\cite{Pas}, \cite{CLb}, One Dimensional Cell decomposition]\label{cdl}
Let $k$ be a number field and $L$ a field over $k$. Let $g_j$ be
polynomials in $1$ variable over $k$. Then  $L((t))$ can be written
as the disjoint union of finitely many $\cL_{DP}(k)$-definable cells
$X_i$ with $\cL_{DP}(k)$-definable presentations $f_i$ and
$\cL_{DP}(k)$-definable centers $c_{i}$ such that for each $i,j$ the
restrictions $\ord g_{j|X_i}$ and $\ac g_{j|X_i}$ factorize through
$f_i$. Moreover, the formulas describing the $X_i$, $f_i$, and $c_i$
can be taken independently of $L$.
\end{lem}
\begin{proof}
This form of Denef - Pas cell decomposition is proven in \cite{CLb}
with $k=\QQ$, where it is a relative cell decomposition relative to
$\ZZ^\ell\times L^\ell$ under the map $(\ord (g_j), \ac (g_j))_j$.
The Lemma for general $k$ follows easily from it by first replacing
the parameters from $k$ by variables, then applying the form of
\cite{CLb} to all the variables, and then plugging in the parameters
for the variables again.
\end{proof}

\begin{cor}\label{lcd}
Let the functions $a_{iz}(K)$, $\alpha_i(K)$, and $\beta_i(K)$ be as
in the proof of Lemma \ref{padicf}. Then there are $M>0$ and a
polynomial $g$ in one variable over $k$ such that for all $a\in k$
with $g(a)\not=0$, all $K\in\cA_{k,M}$ and all $z\in k_K^{n-1}$ one
has that these functions $a_{iz}(K)$, $\alpha(K)$, and $\beta(K)$
are constant functions in the variable $b$ when $b$ runs over
$a+\cM_K$.
\end{cor}

\begin{proof}
Let $g_\ell$ be the polynomials in one valued field variable that
occur in the valued field quantifier free formulas that describe the
$\alpha_i,\beta_i,$ and $A_i$. Apply Lemma \ref{cdl} for the
polynomials $g_\ell$ to find  cells $X_j$ with centers $c_j$ and
presentations $f_j$.
 By quantifier elimination in the language $\cL_{DP}(k)$ and by
examining valued field quantifier free formulas in one free valued
field variables, it follows that the images of the centers $c_j$ are
contained in the zero set of a polynomial $g$ over $k$. Now the
corollary follows from the definition of cells and the evident
ultraproduct argument to go to $K$ in $\cA_{k,M}$ for $M$ big
enough. Namely, for $a\in k$ with $g(a)\not=0$, all $K\in\cA_{k,M}$
for $M$ big, the set $a+\cM_K$ is contained in $X_j(K)$ for some
unique $j$, and $f_j(K)(a+\cM_K)$ is a singleton, where $X_j(K)$ and
$f_j(K)$ are the $K$-interpretations of $X_j$ and $f_j$.
\end{proof}

\section{Estimation of the motivic oscillation index}\label{s:est}

Let $f$ be a polynomial in $n$ variables over a number field $k$.
Let $C_f$, $S_f$, $\delta_f$, and $\varepsilon_f$ be as in the
section on notation. Let $\alpha_f$ be the motivic oscillation index
as before. 

\begin{theorem}\label{fund}
\label{est} Let $f$ be a polynomial in $n$ variables over a number
field $k$. Then
\begin{equation}\label{bounc}
\frac{-n+\delta_f}{2}\leq \alpha_{f}.
\end{equation}
\end{theorem}
\begin{proof}

Since one can always replace $k$ by a finite field extension and $f$
by $f-a$ for some critical value $a\in k$, it is clear that it is
enough to prove that
\begin{equation}\label{bouncbis}
\frac{-n+\varepsilon_f}{2}\leq \alpha_{f}.
\end{equation}
We may suppose that $\varepsilon_f\geq 0$ since otherwise
$\varepsilon_f=-\infty$ and then (\ref{bouncbis}) is trivial. We
then prove the following statement by induction on $n$ when
$\varepsilon_f\geq 0$.
\begin{quote}\textit{
For each $M$ there exist $L$ in $\cA_{k,M}$, a basic step function
$\varphi_L$ on $L^n$ (in the sense of section \ref{sc}), a real
number $s_0$ with
$$\frac{-n+\varepsilon_f}{2}\leq s_0\leq 0,$$
and a complex number $t_0$ with real part $q_L^{-s_0}$ such that
$t_0$ is a nontrivial pole (in the sense of Definition \ref{nont})
of the rational function in $T=q_L^{-s}$ defined for $s>0$ by
$$
Z_{f,L,\varphi_L}(s):=\int_{L^n}\varphi_L(x)|f(x)|^{s}|dx|.
$$
}
\end{quote}
This statement implies (\ref{bouncbis}) (since a nontrivial pole of
$Z_{f,L,\varphi_L}(s)$ corresponds to an effectively appearing
$\lambda$ in (\ref{bounc}) by \cite{DenefVeys}, Proposition 2.7).

The case that $n=1$ is easy by explicit calculation using
factorization of $f$ into irreducible polynomials over the algebraic
closure of $k$.

For $n>1$, first suppose that $\varepsilon_f=0$.
Let $M$ be big enough and choose $L$ in $\cA_{k,M}$ such that
$S_f(k_L)$ is nonempty. By Theorem \ref{tpol} and Remark \ref{tpz}
it follows that $ Z_{f,L,\varphi_{\rm triv}}(s)$ has a nontrivial
pole in $q_L^{-s}$, with $\varphi_{\rm triv}$ the basic step
function which is the characteristic function of $\cO_L^n$. Since
each complex pole $t_0$ of $ Z_{f,L,\varphi_{\rm triv}}(s)$ with
real part $q_L^{-s_0}$ satisfies $-n/2\leq s_0 \leq 0$ when $n>1$ by
results by Segers \cite{Segers}, the case that $\varepsilon_f=0$ and
$n>1$ follows.

Now we treat the case $n>1$ and $\varepsilon_f>0$. There is a
coordinate, say $x_1$, such that for infinitely many points $a$ in
the algebraic closure of $k$, $\varepsilon_{f_a}\geq
\varepsilon_f-1$, with $f_a$ the polynomial in $\hat
x:=(x_2,\ldots,x_n)$ over $k(a)$ defined by $f_a(\hat x)=f(a,\hat
x)$. Indeed, since $\varepsilon_f>0$ there is a coordinate, say
$x_1$, to which an irreducible component of $S_f$ of dimension
$\varepsilon_f$ projects dominantly, and then one compares the
equations defining $S_f$ and the $S_{f_a}$. Let $g$ be a polynomial
in one variable as in Lemma \ref{padicf} for the polynomial $f$.
 Then, by replacing $k$ by a bigger field extension,
we can suppose that there is $c\in k$ such that $g(c)\not=0$ and
such that $\varepsilon_{f_c}\geq \varepsilon_f-1$.
Note that $f_c$ is a polynomial in $n-1$ variables.
 By the induction hypothesis, there exist $L$ in $\cA_{k,M}$, a basic
step function $\psi_L$ on $L^{n-1}$ and a real number $s_0$ with
 $$
\frac{-n+\varepsilon_f}{2}=\frac{-(n-1)+(\varepsilon_f-1)}{2}\leq
\frac{-(n-1)+\varepsilon_{f_c}}{2} \leq s_0\leq 0
 $$
and a complex number with real part $q_L^{-s_0}$ which is a
nontrivial pole of $Z_{f_c,L,\psi_L}(s)$. Now let $\varphi_L$ be the
basic step function on $L^n$ which is the product of the
characteristic function of $c+\cM_L$ with $\psi_L$. Then $t_0$ is
also a nontrivial pole of $ Z_{f,L,\varphi_L}(s)$, since $
Z_{f,L,\varphi_L}(s)$ is a product integral by the above application
of Lemma \ref{padicf} to $f$ and since $g(c)\not=0$. Namely, $
Z_{f,L,\varphi_L}(s)= q_L^{-1}Z_{f_c,L,\psi_L}(s)$. This proves the
above statement and thus the theorem.
 \end{proof}

\begin{remark}
A similar but refined proof than that of Theorem \ref{fund} might
give the stronger result
\begin{equation}\label{bouncd}
 \frac{-n+\delta_f(K)}{2}\leq \alpha_{f,K}
\end{equation}
for all $K\in\cC_{k,M}$ with $M$ sufficiently big and $\delta_f(K)$
the Zariski dimension of the Zariski closure of $C_f(K)$ in
$\AA_k^n$. Note that this implies (\ref{bounc}) since there are
enough $K$ with $\delta_f(K)=\delta_f$.
\end{remark}

\section{Igusa's conjecture modulo $p$ and modulo $p^2$}

\begin{theorem}\label{ti}
Let $f$ be a homogeneous polynomial in $n$ variables over the number
field $k$. Then there exists $M> 0$ and a constant $c$ such that for
each $K$ in $\cC_{k,M}$ and each $y\in K$ with $v_K(y)=-1$ or
$v_K(y)=-2$ one has
\begin{equation}\label{51}
|E_{f,K}(y)|  <  c\, |y|_K^{\alpha_{f}}.
\end{equation}
Moreover, the statement for $y$ with $v_K(y)=-2$ holds for all
polynomials over $k$.
\end{theorem}
\begin{proof}
Let $\delta_f$ be as in the section on notation. First note that for
homogeneous $f$ one has $\delta_f=-\infty$ if and only if $f$ is
linear, and, in the linear case one has $\alpha_f=-\infty$ and
$E_{f,K}(y)=0$ for all $K\in\cC_{k,M}$ for $M$ big enough and all
$y\in K$ with $v_K(y)<0$. The case that $f$ is a constant is
trivial, cf.~Proposition \ref{nontrivialosc}.

Suppose now that the degree of $f$ is $\geq 2$. Thus
$n>\delta_f>-\infty$.
 Let $H_f$ be the (not necessarily reduced)
projectivization in $\PP_k^n$ of the closed subscheme of $\AA^n_k$
given by $f=0$, and let $\delta$ be the (projective) dimension of
the singular locus of the intersection of $H_f$ with the hyperplane
at infinity. Then
\begin{equation}\label{Kaa}
\delta+1=\delta_f.
\end{equation}
Then, by Katz \cite{Katz}, Theorem 4, and by (\ref{Kaa}), there
exist constants $C$ and $M$ such that for all $K$ in $\cC_{k,M}$ and
all $y$ in $K$ with $v_K(y)=-1$ one has
\begin{equation}\label{Ka4}
|E_{f,K}(y)|<C q_K^{\frac{-n+\delta_f}{2}},
\end{equation}
with $q_K$ the number of elements of the residue field of $K$.
Together with (\ref{bounc}) of Theorem \ref{est} this implies the
statement for $y$ with $v_K(y)=-1$.

\par
Now let $f$ be a general polynomial over $k$ in $n$ variables. For
$K\in\cC_{k,M}$ with $M$ big enough, define
$$
A_f(K):=\{x\in\cO_K^n\mid \grad f(x)\equiv 0 \bmod \cM_K\}
$$
and
$$
B_f(K):=\{x\in\cO_K^n\mid \grad f(x)\not\equiv 0 \bmod \cM_K\}.
$$
By 
Noether's Normalization Theorem, the number of elements of the set
$C_f(k_K)$ is $\leq D q_K^{\delta_f}$ for some $D$ independent of
$q_K$ when the characteristic of $k_K$ is big enough. Hence, there
exists some $M$ and $D$ such that for all $K$ in $\cA_{k,M}$ the
measure of $A_f(K)$ against the normalized Haar measure $|dx|$
satisfies
$$
\int_{A_f(K)}|dx|\leq D q_K^{-n+\delta_f}.
$$
 On the other hand, for $M$ big enough, $K$ in $\cA_{k,M}$, and $y\in K$ with $v_K(y)=-2$,
$$
0=\int_{B_f(K)}\psi_K(yf(x))|dx|_K,
$$
by Hensel's Lemma and the fact that the sum
$$\sum_{z\in\FF_q}\psi_q(z)$$ equals zero for any nontrivial additive
character $\psi_q$ on $\FF_q$.
 Hence, for such $M$ and $K$, and for all $y$ in $K$ with $v_K(y)=-2$,
\begin{eqnarray} \nonumber \abs{E_{f,K}(y) } & =  & \abs{ \int_{A_f(K)}\psi_K(yf(x))|dx|_K
 + \int_{B_f(K)}\psi_K(yf(x))|dx|_K }\\
\nonumber  & = & \abs{ \int_{A_f(K)}\psi_K(yf(x))|dx|_K}\\
\nonumber   & \leq  & \int_{A_f(K)}|dx|_K\\
\nonumber   & \leq  & D q_K^{-n+\delta_f}\\
    & \leq  & D |y|^{\alpha_f},\label{412}
\end{eqnarray}
where $D$ is independent of $K$ and where inequality (\ref{412})
follows from $|y|=q_K^2$ and (\ref{bounc}).
\end{proof}

\begin{cor}\label{ci}
Igusa's Conjecture \ref{igusacon} with argument of order $-1$ or
$-2$ holds for all homogeneous polynomials $f$ over $k$. Namely,
for $f$ a homogeneous polynomial in $n$ variables over $k$ and
$\pi=\pi_f$ a resolution as in the introduction, there exists a
constant $c$ such that for each $p$-adic completion $K$ of $k$ and
each $y\in K$ with $v_K(y)=-1$ or $v_K(y)=-2$ one has
$$
|E_{f,K}(y)|  <  c\, |y|_K^{\alpha(\pi)}.
$$
\end{cor}
\begin{proof}
This follows from Theorem \ref{ti} and (\ref{lin}) for big residue
characteristics, and from (\ref{a2}) for small residue
characteristics.
\end{proof}

\begin{remark}
The factor $|v_K(y)|^{n-1}$ featuring in the original conjecture
clearly becomes unnecessary when one focuses on $y$ with
$v_K(y)=-1$ and $v_K(y)=-2$. When $v_K(y)$ varies over $\ZZ$, the
factor $|v_K(y)|^{n-1}$ can in general not be omitted.
\end{remark}

\section{Nonhomogeneous polynomials}

Let $f$ be a polynomial in $n$ variables over a number field $k$.
Let $H_f$ be the (not necessarily reduced) projectivization in
$\PP_k^n$ of the scheme $f=0$ in $\AA^n_k$. Let $X_f$ be the
scheme-theoretic intersection of $H_f$ with the hyperplane at
infinity. Let $\delta$ be the dimension of the singular locus of
$X_f$ when this singular locus is nonempty, and put $\delta=-1$ when
$X_f$ is smooth. We say that $f$ satisfies Tameness Condition
(\ref{cond}) if
 \begin{equation}\label{cond}
\delta_f \geq \delta+1.
 \end{equation}

Note that by (\ref{Kaa}) a homogeneous polynomial of degree $\geq 2$
always satisfies Tameness Condition (\ref{cond}). Under this
Tameness Condition we get the analogue of Theorem \ref{ti}:

\begin{theorem}\label{gig}
Let $f$ be a polynomial in $n$ variables over a number field $k$.
Suppose that $f$ satisfies Tameness Condition (\ref{cond}). Then
there exists a number $M$ and a constant $c$ such that for each
$K\in\cC_{k,M}$ and each $y\in K$ with $v_K(y)=-1$ one has
\begin{equation}\label{61}
|E_{f,K}(y)|< c\, |y|^{\alpha_{f}},
\end{equation}
\end{theorem}
\begin{proof}
Same proof as for Theorem \ref{ti}, using $\delta+1$ instead of
$\delta_f$ in (\ref{Ka4}), which holds by the same Theorem 4 of
\cite{Katz}.
\end{proof}

\begin{example}\label{exa}
For $f(x_1,x_2)=x_1^2x_2-x_1$, one has
$\alpha(\pi)=\alpha_f=-\infty$, but $E_{f,K}(y)\not = 0$ for any
$K\in\cC_{\QQ,2}$  when $v_K(y)=-1$. Hence, a literal analogue of
Conjecture \ref{igusacon} would not make sense for general
polynomials. Note that $f$ does not satisfy Tameness Condition
(\ref{cond}).
\end{example}

\section{Igusa's conjecture revisited}\label{sglobal}

We come back to the sums $E_f(N)$ from the introduction and study
the growth of $|E_f(N)|$ with $N$. As noted before, it is sufficient
to bound $|E_f(p^m)|$ in terms of primes $p$ and positive integers
$m$. We focus on $m>1$ since the case $m=1$ is rather well
understood by work by Deligne, Laumon, Katz, and others.

In general one has the following situation which follows from
\cite{DenefBour} or \cite{CLexp}. Let $f$ be any polynomial over
$\QQ$ in $n$ variables. Then there exist constants $\beta,\gamma\geq
0$ such that for all big enough primes $p$ and all $m\geq \gamma$,
one has
\begin{equation}\label{beta}
|E_f(p^m)| < p^\beta m^{n-1}  p^{\alpha_f m},
\end{equation}
with $\alpha_f$ the motivic oscillation index.

A natural sharpening of Igusa's Conjecture \ref{igusacon} is the
conjecture that one can take $\beta=\gamma=0$ in (\ref{beta}) when
$f$ is homogeneous. (Conjecture \ref{igusacon} follows from this
sharpening by (\ref{a2}) and by Corollary \ref{cor1}.) An
alternative sharpening of Igusa's conjecture by Denef and Sperber in
\cite{DenSper} is based on complex oscillation indices and might be
equivalent with our sharpening (that is, the supremum of the complex
oscillation indices of $f$ around each of the points of $\CC^n$
might equal $\alpha_f$).

We introduce the following mysterious invariant of $f$ which might
play a role in understanding $E_f(p^m)$ for general $f$:

\begin{definition}
Let $f$ be a polynomial in $n$ variables over a number field $k$.
Define $\beta_f$ as the limit over $M$ and $\gamma$ of the infima of
all real $\beta\geq 0$ such that
$$
|E_{f,K}(y)|< q_K^{\beta}|v_K(y)|^{n-1}|y|_K^{\alpha_f}
$$
for all $y\in K$ with $v_K(y)<-\gamma$ and all $K\in\cA_{k,M}$. Call
$\beta_f$ the  \emph{flaw of $f$ w.r.t.~$\alpha_f$}.
\end{definition}

For non-homogeneous $f$, $\beta_f$ is not even conjecturally
understood.

\subsection*{Acknowledgment}
I would like to thank J.~Denef, F.~Loeser, W.~Veys, D.~Segers, and
J.~Nicaise for fruitful discussions during the preparation of this
paper. I am indebted to J.~Denef for his suggestions on the
definition of the motivic oscillation index.


\bibliographystyle{amsplain}
\bibliography{anbib}

\end{document}